\documentclass[a4paper,10pt]{article}

\usepackage{amssymb, amsmath, latexsym, amsthm, cite}

\newcommand{\Q}{\mathbb{Q}}

\newcommand{\F}{\mathbb{F}}

\renewcommand{\b}[1]{{\bf #1}}

\newcommand{\x}{\b{x}}

\newcommand{\y}{\b{y}}

\newcommand{\z}{\b{z}}

\newtheorem{theorem}{Theorem}

\newtheorem{corollary}{Corollary}

\newtheorem{lemma}{Lemma}

\begin{document}

\title{Simultaneous zeros of a Cubic  \\ and Quadratic form}

\author{Jahan Zahid}

\date{}
\maketitle

\section{Introduction} 
\label{intro}

Consider a system of forms
\[ \b{F}(\x) = (F_{1}(\x),\ldots,F_{r}(\x)) \]
of degrees $d_{1},\ldots,d_{r}$ respectively in the variables $\x = (x_{1},\ldots,x_{n})$ over a $\mathfrak{p}$-adic field $K$. It had been conjectured by Artin \cite[Preface]{Art65} that $\b{F}$ necessarily has a non-trivial zero in $K$ provided $n > \sum d_{i}^{2}$. It should be noted that there exist systems with $n = \sum d_{i}^{2}$ which have only the trivial $\mathfrak{p}$-adic zero, so this is the best we can hope for. 

Artin's conjecture has been verified in only a handful of cases. For example it is a classical result due to Hasse \cite{Has24} that every quadratic form with $n > 4$ variables has a non-trivial $\mathfrak{p}$-adic zero. The case of cubic forms was settled independently by Dem'yanov \cite{Dem50}, Lewis \cite{Lew52} and Springer \cite{Spr55}. Dem'yanov \cite{Dem56} and later Birch, Lewis \& Murphy \cite{BirLewMur62} proved the conjecture for a system of two quadratic forms. The purpose of this paper is to prove the next case of the conjecture for two forms, provided we have a large enough residue class field. More precisely we shall prove

\begin{theorem} \label{mainres}
Any system of a cubic and quadratic form in at least $14$ variables defined over $K$, has a non-trivial zero in $K$ provided the cardinality of the residue class field exceeds $293$.
\end{theorem}  

It should be noted that Artin's conjecture was shown to be false in general by Terjanian \cite{Ter66}, who found a counterexample of a quartic form in 18 variables with no zero in $\Q_{2}$. If  however $|K : \Q_{p}| = e$ and $\b{d}=(d_{1},\ldots,d_{r})$, then by a remarkable theorem of Ax \& Kochen \cite{AxKoc65}, there exists an integer $p(\b{d},e)$ such that any system $\b{F}$ with $n > \sum d_{i}^{2}$ has a non-trivial zero provided the characteristic of the residue class field exceeds $p(\b{d},e)$. 

We remark that Theorem \ref{mainres} is stronger than anything we can deduce from the Ax--Kochen theorem for a number of reasons. Firstly we have an explicit bound on the cardinality of the residue class field for which Artin's conjecture is true. Secondly we have a condition depending on the cardinality of the residue class field, rather than the characteristic. Consequently we are now able to say that Artin's conjecture holds for a cubic and quadratic form over any unramified extension of $\Q_{p}$ of degree at least $9$. Where it was not possible to make this deduction before. 

As an outline to prove Theorem \ref{mainres} we shall generalise a $\mathfrak{p}$-adic minimization procedure due to Schmidt \cite{Sch80} to hold for systems of forms of arbitrary degrees. We shall then derive some Geometric information of the system over the residue class field, for those systems which terminate in the minimization process. This will allow us to find a non-singular zero in the residue class field to which we can apply Hensel's Lemma.
\vspace{.1in}

\b{Acknowledgments:} This work forms part of the authors doctoral thesis at the University of Oxford. I very gratefully acknowledge the financial support I received from EPSRC. 

I would also like to thank my supervisor Prof. Roger Heath-Brown, for suggesting this problem and his excellent guidance over the last few years. I have also benefited from numerous conversations with Dr. Damiano Testa and Prof. Trevor Wooley, to whom I am glad to express my gratitude.

\section{Some preliminaries} 

Let $\mathcal{O}_{K}$ denote the ring of integers of $K$, and denote the residue class field by $\F_{q}$. Let $\pi$ denote a uniformizer for $\mathcal{O}_{K}$. If $\alpha \in K-\{0\}$, we may write $\alpha = \pi^{s}u$, where $u$ is a unit in $\mathcal{O}_{K}$. We define the $\pi$-adic order $v(\cdot)$ by setting $v(\alpha) = s$. We also define the $\pi$-adic valuation $| \cdot |$ by setting $| \alpha | = p^{-s}$, where $p$ denotes the characteristic of the residue class field $\F_{q}$.
Recall that 
\[ F_{1}(\x),\ldots,F_{r}(\x) \in K[\x] \]
denotes an arbitrary system of forms of degrees $d_{1},\ldots,d_{r}$ in the variables $\x = (x_{1},\ldots,x_{n})$ over $K$. We shall assume that $d_{1} \ge \ldots \ge d_{r}$ unless we state otherwise and for brevity write the above system of forms as $\b{F}$. We are interested in determining the existence of a point $\x \in K^{n} - \{ \b{0} \}$ such that $\b{F}(\x) = \b{0}$. Clearly we may assume that the coefficients of the forms $\b{F}$ and the variables $\x$ are in $\mathcal{O}_{K}$, since this does not affect existence of a zero. 

By a slight abuse of notation we write $GL(n,\mathcal{O}_{K})$ to denote the set of $(n \times n)$-matrices over $\mathcal{O}_{K}$ with non-zero (rather than unitary) determinant. So let $\tau \in GL(n,\mathcal{O}_{K})$ and write $\b{F}_{\tau}$ to denote $\b{F}(\tau \x) = (F_{1}(\tau \x),\ldots,F_{r}(\tau \x))^{t}$. We also write $T = (t_{ij})$ to denote the $(r \times r)$ upper triangular matrix with entries 
\begin{eqnarray} \label{uptrimat}
t_{ij}(\x) = \pi^{-c_{i}}G_{ij}(\x), 
\end{eqnarray}
where $c_{i} \ge 0$ and $G_{ij} \in \mathcal{O}_{K}[\x]$ to denote an arbitrary form of degree 
\[ \deg F_{i} - \deg F_{j} \ge 0, \] 
for $1 \le i < j \le r$ and define the diagonal terms $G_{ii} := 1$. 

Let $\Omega= (\omega_{1},\ldots,\omega_{r})$ where $\omega_{i}>0$ for each $i$. Then we write 
\[ \b{F} \mathop{\succ}_{\Omega} \b{F'}, \]
if $\b{F}$ and $\b{F'}$ are both defined over $\mathcal{O}_{K}$ and
\[ \b{F'} = T\b{F}_{\tau} \]
with
\begin{eqnarray} \label{bot-cond} 
\sum c_{i} \omega_{i} - s > 0 
\end{eqnarray}
where $s = v(\det \tau)$ and the $c_{i}$ are as in (\ref{uptrimat}). If $\b{F'} = T\b{F}_{\tau}$, it is clear that if $\b{F'}$ has a zero if and only  if $\b{F}$ has a zero. 

We say that $\b{F}$ is $\Omega$-bottomless if there is an infinite chain
\[ \b{F} \mathop{\succ}_{\Omega} \b{F}^{(1)} \mathop{\succ}_{\Omega} \b{F}^{(2)} \mathop{\succ}_{\Omega} \cdots \]
otherwise, $\b{F}$ will be called $\Omega$-bottomed. We also say that $\b{F}$ is $\Omega$-reduced if there does not exist any $\b{F'}$ such that
\[ \b{F} \mathop{\succ}_{\Omega} \b{F'}. \]

We say that two systems $\b{F}$ and $\b{F'}$ are equivalent, if both systems are defined over $\mathcal{O}_{K}$ and 
\[ \b{F'} = T\b{F}_{\tau} \]
where $c_{i}=0$ for all $1 \le i \le r$ and $v(\det \tau) = 0$ in $T$ and $\tau$ as in (\ref{uptrimat}). The order $o(\b{F})$, of a system $\b{F}$ is the least positive integer $m$ such that $\b{F}$ is equivalent to a system that contains $m$ variables explicitly. We also define the $h$-invariant for a system $\b{F}$, denoted $h(\b{F})$ as the least integer $h$ such that we can write
\begin{equation} \label{h-inv}
F_{i}(\x) = x_{1}H_{i1}(\x) + \ldots + x_{h}H_{ih}(\x)  \pmod \pi, 
\end{equation}
for all $1 \le i \le r$ and all systems equivalent to $\b{F}$. Note that since the $F_{i}$ are defined over $\mathcal{O}_{K}$, considering them modulo $\pi$ is well defined.

Given any set $S = \{ e_{1},\ldots,e_{s} \}$ of positive integers we define $v_{S}$ as the least integer $v$ such that every system consisting of $s$ forms of degrees $e_{1},\ldots,e_{s}$ have a non-trivial zero provided the number of variables in the system is at least $v_{S}$. If $\phi$ denotes the empty set we define $v_{\phi} := 1$. We remark here that it is due to a classical theorem of Brauer \cite{Bra45}, the number $v_{S}$ is always finite. 

Although the next theorem is likely to have further applications, it will for the purpose of this paper play a crucial part in the minimisation procedure for a system of a cubic and quadratic form.

\begin{theorem} \label{b-less}
Let $S \subset \{ d_{1},\ldots,d_{r} \}$ denote any subset of cardinality $r-1$ with indexing set $I$ such that $v_{S}$ is maximal. Let $j \not\in I$ then provided 
\begin{eqnarray} \label{botass}
n \ge v_{S}+d_{j}^{2} 
\end{eqnarray}
there exists some $\Omega = (\omega_{1},\ldots,\omega_{r})$ such that $\omega_{i} > d_{i}$ for each $1 \le i \le r$ and such that every $\Omega$-bottomless system $\b{F}$ defined over $\mathcal{O}_{K}$ has a non-trivial $\mathfrak{p}$-adic zero.
\end{theorem}


\section{Proof of Theorem \ref{b-less}}

Since the field $K$ has characteristic $0$, given any form $F$ of degree $d$ there is a unique form $M_{F}(\x_{1},\ldots,\x_{d})$ which is linear in each vector $\x_{j}$ and which is symmetric in $\x_{1},\ldots,\x_{d}$, such that

\[ F(\x) = M_{F}(\x,\ldots,\x). \]

Let $\b{e}_{1}=(1,0,\ldots,0), \b{e}_{2}, \ldots, \b{e}_{n}$ be unit vectors. We say that $\b{F} = (F_{1},\ldots,F_{r})$ is $\Omega$-special if there are non-negative integers $a_{1},\ldots,a_{n}$ and $b_{1},\ldots,b_{r}$ with
\begin{eqnarray} \label{spec1}
a_{1}+\cdots+a_{n} < \omega_{1} b_{1} + \cdots + \omega_{r} b_{r}
\end{eqnarray}
such that 
\begin{eqnarray} \label{spec2}
M_{F_{i}}(\b{e}_{j_{1}},\ldots,\b{e}_{j_{d_{i}}}) = 0
\end{eqnarray}
for each $1 \le i \le r$ and $d_{i}$-tuple $(j_{1},\ldots,j_{d_{i}})$ for which, 
\begin{eqnarray} \label{spec3}
a_{j_{1}}+\cdots+a_{j_{d_{i}}}<b_{i}. 
\end{eqnarray}

Note the following important correspondence between $\Omega$-bottomless systems and $\Omega$-special systems.

\begin{theorem} \label{bot-spec}
Every $\Omega$-bottomless system is equivalent to a $\Omega$-special system.
\end{theorem}
We prove this in due course. Now we make use of this to prove Theorem \ref{b-less}.

\begin{proof} [Proof of Theorem ]

By Theorem \ref{bot-spec} we may suppose that $\b{F}$ is $\Omega$-special. For ease of notation we may assume that 
\begin{eqnarray} \label{blessord}
a_{1} \le \ldots \le a_{n} \quad \mbox{and} \quad \frac{b_{1}}{d_{1}} \le \ldots \le \frac{b_{r}}{d_{r}}, 
\end{eqnarray}
dropping any previous ordering we had on $d_{1},\ldots,d_{r}$. If there is a subset $S \subset \{ d_{1},\ldots,d_{r} \}$ with indexing set $I$ such that  
\[ d_{i}a_{v_{S}} < b_{i}  \quad \mbox{for all} \; \; i \not\in I \]
then $\b{F}$ has a non-trivial zero. For if such a subset $S$ exists then by (\ref{spec3}) one has
\[ M_{F_{i}}(\b{e}_{j_{1}},\ldots,\b{e}_{j_{d_{i}}}) = 0 \quad \mbox{for all} \; \; i \not\in I \]
and every $1 \le j_{1},\ldots,j_{d_{i}} \le v_{S}$. Therefore the system $(F_{i})_{i \not\in I}$ vanishes on the $v_{S}$-dimensional subspace spanned by $\{ \b{e}_{1},\ldots,\b{e}_{v_{S}} \}$ and on this subspace we can find a zero of $(F_{i})_{i \in I}$ and therefore a zero of $\b{F}$. Consequently we may assume that for each $S \subset \{ d_{1},\ldots, d_{r} \}$ there exists some $i$, which by the ordering (\ref{blessord}) we may assume to be $\min \{1 \le i \le d : i \not\in I \}$ such that 
\[ d_{i}a_{v_{S}} \ge b_{i}. \]
We define $S_{0}:= \phi$ and $S_{i} := \{ d_{1},\ldots,d_{i} \}$ for $1 \le i \le r-1$ and for ease of notation write $w_{i} = v_{S_{i-1}}$ for $1 \le i \le r$. Then it follows that
\begin{eqnarray} \label{daw}
d_{i}a_{w_{i}} \ge b_{i} 
\end{eqnarray}
for every $1 \le i \le r$. Note that by assumption (\ref{botass}) we have that 
\[ n \ge w_{r} + d_{r}^{2} \ge d_{1}^{2}+ \cdots + d_{r}^{2}+1. \] 
Moreover for each $0 \le i \le r-1$ we claim that
\begin{eqnarray} \label{wd}
w_{i+1} - w_{i} \ge d_{i}^{2}. 
\end{eqnarray}
For if we let $\b{T}=(F_{1},\ldots,F_{i-1})$ denote a system with $w_{i}-1$ variables with only the trivial zero and $F_{i}$ denote a form in $d_{i}^{2}$ variables with only the trivial zero with its variables distinct from the variables in $\b{T}$ then it is clear that the system $\b{T} \cup (F_{i})$ has only the trivial zero. Therefore it follows that $w_{i+1} \ge w_{i} + d_{i}^{2}$ as claimed. By  (\ref{blessord}), (\ref{daw}) and (\ref{wd}) it follows that
\begin{eqnarray*} 
a_{1} + \cdots + a_{n} & \ge & d_{1}^{2}a_{w_{1}} + d_{2}^{2}a_{w_{2}}+ \cdots + d_{r}^{2}a_{w_{r}} + a_{n} \\
& \ge & \Big(d_{1}^{2}+\frac{1}{r} \Big)a_{w_{1}} + \Big(d_{2}^{2}+\frac{1}{r} \Big)a_{w_{2}}+ \cdots + \Big(d_{r}^{2}+\frac{1}{r} \Big)a_{w_{r}} \\
& \ge & (d_{1}+\epsilon)b_{1} + (d_{2}+\epsilon)b_{2} + \cdots + (d_{r}+\epsilon)b_{r}
\end{eqnarray*}
taking $\epsilon = (rd_{max})^{-1}$, where $d_{max}:= \max \{d_{1},\ldots,d_{r} \}$. Hence if we let $\omega_{i} = d_{i}+\epsilon$ for $1 \le i \le r$, it follows that every $\Omega$-special system must have a non-trivial zero completing the proof of the theorem.
\end{proof}

We shall now proceed by proving Theorem \ref{bot-spec}, generalising where appropriate the method of Schmidt \cite{Sch80}. Given two systems $\b{F}$ and $\b{F'}$ defined over  $\mathcal{O}_{K}$, we write
\begin{eqnarray} \label{high}
\b{F} \mathop{\succ}_{\Omega}^{k} \b{F'} 
\end{eqnarray}
if $k \ge 1$ and 
\[ \b{F} \mathop{\succ}_{\Omega} \b{F'} \quad \; \;  \]
with the condition (\ref{bot-cond}) strengthened to 
\[ \sum c_{i} \omega_{i} - (s+k) \ge 0. \]
The system $\b{F}$ shall be called $\Omega$-high if for every $k$ there is a $\b{F'}$ such that (\ref{high}) holds. Note the following lemma. 

\begin{lemma} \label{bot-high}
Suppose $\b{F}$ is a $\Omega$-bottomless system, then it is $\Omega$-high.
\end{lemma}

\begin{proof}
Fix a $k \ge 1$, then if $\b{F}$ is $\Omega$-bottomless there exists an infinite chain
\[ \b{F} = \b{F}^{(1)} \mathop{\succ}_{\Omega} \b{F}^{(2)} \mathop{\succ}_{\Omega} \ldots \mathop{\succ}_{\Omega} \b{F}^{(k)} \mathop{\succ}_{\Omega} \cdots. \]
For each $m \ge 1$ we may write 
\[ \b{F}^{(m+1)} = T_{m}\b{F}_{\tau_{m}}^{(m)} \]
where $\tau_{m} \in GL(n,\mathcal{O}_{K})$ with $v(\det \tau_{m}) = s_{m}$ and $T_{m} = (t_{ij,m})$ is an $(r \times r)$ upper triangular matrix with entries 
\begin{eqnarray*} 
t_{ij,m}(\x) = \pi^{-c_{i,m}}G_{ij,m}(\x), 
\end{eqnarray*}
for $1 \le i \le j \le r$ where $c_{i,m} \ge 0$ and $G_{ij,m} \in \mathcal{O}_{K}[\x]$ is a form of degree\footnote{we assume once again that $d_{1} \ge d_{2} \ge \ldots \ge d_{r}$} 
\[ \deg F_{i} - \deg F_{j} \ge 0 \] 
with $G_{ii}=1$. Let $Q$ be any positive integer such that $\omega_{i} \ge \frac{1}{Q}$ for every $i$, then by condition (\ref{bot-cond}) we have that
\[ \sum_{i} c_{i,m} \omega_{i} - (s_{m} +\tfrac{1}{Q}) \ge 0 \]
for each $1 \le m \le r$. If $T$ denotes the product of $T_{m}$ and $\tau$ denotes the product of $\tau_{m}$ for $1 \le m \le kQ$, then we have that $\b{F}^{(kQ+1)} = T \b{F_{\tau}}$. Crucially we also have that 
\[ \sum_{i} \Big( \sum_{m} c_{i,m} \Big) \omega_{i} - \Big( \sum_{m} s_{m} +k \Big) \ge 0. \]
Therefore 
\[ \b{F} \mathop{\succ}_{\Omega}^{k} \b{F}^{(kQ+1)}  \] 
as required.
\end{proof}

We shall now note two Lemmata, the proofs of which can be found in Schmidt's paper \cite[Lemmata 8 and 10]{Sch80}.

\begin{lemma} \label{lemma8}
Let $A_{1},\ldots,A_{l}$ and $B_{1},\ldots,B_{m}$ be linear forms with integer coefficients in the vector $\x = (x_{1},\ldots,x_{n})$. Let $\x_{1},\x_{2},\ldots$ be a sequence of vectors with 
\[ A_{i}(\x_{k}) \ge 0 \quad (i=1,\ldots,l; \; k=1,2,\ldots). \]
Then there exists a subsequence, $\y_{1},\y_{2},\ldots$ say, a constant $B$ and an integer vector $\b{a}$ with
\[ A_{i}(\b{a}) \ge 0 \quad (i=1,\ldots,l) \]
such that
\[ \lim_{k \to \infty} B_{j}(\y_{k}) = +\infty \quad \mbox{for $j$ with} \; B_{j}(\b{a})>0 \] 
and
\[ B_{j}(\y_{k}) \le B \quad \mbox{for $j$ with} \; B_{j}(\b{a}) \le 0. \]
\end{lemma}

Before stating the next lemma we need to introduce some terminology. If $\b{a}_{1},\ldots,\b{a}_{s}$ are linearly independent vectors in $K^{n}$, we call the set of linear combinations
\[ c_{1}\b{a}_{1}+\cdots+c_{s}\b{a}_{s} \]
a lattice, where $c_{i} \in \mathcal{O}_{K}$. We call $\b{a}_{1},\ldots,\b{a}_{s}$, a basis of the lattice. 

\begin{lemma} \label{lemma10}
Suppose $M$ is a sublattice of $\Lambda$. Then there exists a basis $\b{l}_{1},\ldots,\b{l}_{s}$ of $\Lambda$ and a basis of $\b{m}_{1},\ldots,\b{m}_{s}$ of $M$ such that 
\[ \b{m}_{1} = \pi^{a_{1}}\b{l}_{1},\ldots,\b{m}_{s} = \pi^{a_{s}}\b{l}_{s} \]
for some non-negative integers $a_{1},\ldots,a_{s}$.
\end{lemma}

\begin{proof} [Proof of Theorem $\ref{bot-spec}$]
Suppose $\b{F}$ is $\Omega$-bottomless, then by Lemma \ref{bot-high} it is $\Omega$-high. Hence for every $k \ge 1$ there are maps $T_{k}$ and $\tau_{k} \in GL(n,\mathcal{O}_{K})$ such that $T_{k}\b{F}_{\tau_{k}}$ is defined over $\mathcal{O}_{K}$ and 
\begin{equation} \label{highcond} 
\sum_{i} c_{i,k} \omega_{i} - (s_{k} +k) \ge 0
\end{equation}
where $v(\det \tau_{k}) = s_{k}$ and $T_{k} = (t_{ij,k})$ is an $(r \times r)$ upper triangular matrix with entries 
\begin{eqnarray*} 
t_{ij,k}(\x) = \pi^{-c_{i,k}}G_{ij,k}(\x), 
\end{eqnarray*}
for $1 \le i \le j \le r$ where  $c_{i,k} \ge 0$ and $G_{ij,k} \in \mathcal{O}_{K}[\x]$ is a form of degree 
\[ \deg F_{i} - \deg F_{j} \ge 0 \] 
with $G_{ii}=1$. For ease of notation we write the $i$th row of the vector
\[ diag (\pi^{c_{1,k}},\ldots,\pi^{c_{r,k}}) T_{k} \b{F} (\x) \]
as
\[ R_{i,k}(\x) := \sum_{j=i}^{r} G_{ij,k}(\x)F_{j}(\x) \]
for $1 \le i \le r$. If $\Lambda_{k}$ denotes the lattice $\tau_{k} \mathcal{O}_{K}^{n}$ then by assumption
\[ \pi^{-c_{i,k}} R_{i,k} (\x) \in \mathcal{O}_{K} \]
for all $1 \le i \le r$ and every $\x \in \Lambda_{k}$. By Lemma \ref{lemma10}, $\Lambda_{k}$ has a basis 
\begin{equation} \label{lambas}
\Lambda_{k} : \pi^{u_{1}}\b{u}_{1},\ldots, \pi^{u_{n}}\b{u}_{n} 
\end{equation}
where $\b{u}_{1},\ldots,\b{u}_{n}$ is a basis of $\mathcal{O}_{K}^{n}$. Next we let $M_{R_{i,k}}$ denote the multilinear forms associated with $R_{i,k}$ for $1 \le i \le r$, then there exists some fixed non-negative integer $\gamma$ such that 
\[ M_{R_{i,k}} (\x_{1},\ldots,\x_{d_{i}}) \in \pi^{-\gamma} \mathcal{O}_{K} \]
for all $1 \le i \le r$ and any $\x_{1},\ldots,\x_{d_{i}} \in \Lambda_{k}$. Consequently taking the basis vectors of $\Lambda_{k}$ (\ref{lambas}) we get
\[ \pi^{-c_{i,k}} M_{R_{i,k}} (\pi^{u_{j_{1}}}\b{u}_{j_{1}},\ldots,\pi^{u_{j_{d_{i}}}}\b{u}_{j_{d_{i}}}) \in \pi^{-\gamma} \mathcal{O}_{K} \]
or 
\[ |M_{R_{i,k}} (\b{u}_{j_{1}},\ldots,\b{u}_{j_{d_{i}}})| \le p^{\gamma - (c_{i,k} - u_{j_{1}} - \ldots - u_{j_{d_{i}}})} \]
for all $1 \le j_{1},\ldots,j_{d_{i}} \le n$. 

Note that the $u_{i},\b{u}_{i},c_{i,k}$ and $R_{i,k}$ all depend on $k$. Also note that since $\b{u}_{1},\ldots,\b{u}_{n}$ is a basis of $\mathcal{O}_{K}^{n}$ we must have that $| \det (\b{u}_{1},\ldots,\b{u}_{n}) | = 1$. By the compactness of $\mathcal{O}_{K}^{n}$ there must exist a subsequence of the sequence of integers $k=1,2,\ldots,$ such that on this subsequence $\b{u}_{1},\ldots,\b{u}_{n}$ tend respectively to $\b{a}_{1},\ldots,\b{a}_{n}$ and the forms $R_{1,k},\ldots,R_{r,k}$ tend respectively to the forms $R_{1},\ldots,R_{r} \in \mathcal{O}_{K}(\x)$. To be clear the system $(R_{1},\ldots,R_{r})$ has a zero if and only if the system $\b{F}=(F_{1},\ldots,F_{r})$ has a zero, since $G_{ii}=1$ for each $1 \le i \le r$ in $R_{i,k}$. Moreover we have that $|\det (\b{a}_{1},\ldots,\b{a}_{n})| = 1$, hence $\b{a}_{1},\ldots,\b{a}_{n}$ is a basis of $\mathcal{O}_{K}^{n}$. There exists a map $\sigma$ defined over $\mathcal{O}_{K}$ such that 
\[ \sigma \b{e}_{i} = \b{a}_{i}, \quad \mbox{for $1 \le i \le n$}. \]
With each $k$ in our subsequence we define the vector 
\[ (A_{1},\ldots,A_{n+r})  = (u_{1},\ldots,u_{n},c_{1,k},\ldots,c_{r,k}) \]
We also define
\[ B = \sum_{i=1}^{r} c_{i,k} \omega_{i} - \sum_{j=1}^{n} u_{j} \]
and 
\[ B_{i} = c_{i,k} - (u_{j_{1}} + \cdots + u_{j_{d_{i}}}) \]
for $1 \le j_{1},\ldots,j_{d_{i}} \le n$ and $1 \le i \le r$.  We apply Lemma \ref{lemma8} to the forms $A_{i},B,B_{j}$, but before we do this note that the form $B$ tends to $+ \infty$. This is because $s_{k} = u_{1}+\cdots+u_{n}$ and so by equation (\ref{highcond}) one has
\[  \sum_{i} c_{i,k} \omega_{i} - \sum_{j} u_{j} \ge k. \]
Now by Lemma \ref{lemma8} there is a vector
\[ \b{a} = (a_{1},\ldots,a_{n},b_{1},\ldots,b_{r}) \]
with non-negative integer components, such that 
\[ \omega_{1} b_{1} + \cdots + \omega_{r} b_{r} - (a_{1} + \cdots + a_{n}) > 0. \]
Moreover we have that $B_{i}$ tends to $+\infty$ for all values $i$ and $j_{1},\ldots,j_{d_{i}}$ for which
\begin{equation} \label{Mcond} 
b_{i} - (a_{j_{1}} + \cdots + a_{j_{d_{i}}}) > 0. 
\end{equation}
Taking the limit we obtain
\[ M_{R_{i}} (\b{a}_{j_{1}},\ldots,\b{a}_{j_{d_{i}}}) = 0 \]
for all $j_{1},\ldots,j_{d_{i}}$ satisfying (\ref{Mcond}). If we set $\b{G} = \b{R}_{\sigma}$ where $\b{R} = (R_{1},\ldots,R_{r})$, then $\b{G}$ satisfies precisely the conditions required to be $\Omega$-special. Finally since $\b{R}$ is equivalent to $\b{F}$, we deduce that $\b{F}$ is equivalent to a $\Omega$-special system as required.
\end{proof}


\section{Preliminaries for a Cubic and Quadratic form}

Throughout this section and subsequent sections $\b{F}=(F,G)$ will denote a system of a cubic and quadratic form in $n$ variables defined over the ring of integers $\mathcal{O}_{K}$ of some $\mathfrak{p}$-adic field $K$. Let $\pi$ denote a uniformizer for $\mathcal{O}_{K}$ and let $\F_{q}$ be the residue class field of $K$, where $q$ denotes its cardinality.

Note the following corollary of Theorem \ref{b-less}.

\begin{corollary} 
Let $(F,G)$ denote an arbitrary system of a cubic and quadratic form in $n \ge 14$ variables. Then there exists some $\omega_{1}>3$ and $\omega_{2}>2$ such that every $(\omega_{1},\omega_{2})$-bottomless system $(F,G)$ has a non-trivial $\mathfrak{p}$-adic zero.
\end{corollary}

Therefore it follows that to prove Theorem \ref{mainres} it is sufficient to consider $(\omega_{1},\omega_{2})$-reduced systems for some $\omega_{1}>3$ and $\omega_{2}>2$. 

Recall that we say that two systems $\b{F}$ and $\b{F'}$ are equivalent, if both systems are defined over $\mathcal{O}_{K}$ and 
\[ \b{F'} = T\b{F}_{\tau} \]
where 
\begin{displaymath}
T =  \left( \begin{array}{cc}
1 & L \\
0 & 1 
\end{array} \right),
\end{displaymath} \\
for some linear form $L \in \mathcal{O}_{K}[\x]$ and $\tau \in GL(n,\mathcal{O}_{K})$ with $v(\det \tau) = 0$. The order $o(\b{F})$, of a system $\b{F}$ is the least positive integer $m$ such that $\b{F}$ is equivalent to a system that contains $m$ variables explicitly. Also recall the $h$-invariant of $\b{F}$, denoted $h(\b{F})$ is the least integer $h$ such that we can write
\[ F(\x) = x_{1}H_{11}(\x) + \ldots + x_{h}H_{1h}(\x)  \pmod \pi, \]
and 
\[ G(\x) = x_{1}H_{21}(\x) + \ldots + x_{h}H_{2h}(\x)  \pmod \pi, \]
for all systems equivalent to $\b{F}$. We similarly define the $h$-invariant of a single form in the obvious way. Note here that since $\b{F}$ is defined over $\mathcal{O}_{K}$, considering it modulo $\pi$ is well defined. Note the following lemma which will play a crucial part in our proof.

\begin{lemma} \label{reduced}
Suppose $\b{F} = (F,G)$ is $(\alpha,\beta)$-reduced, for some $\alpha > 3$ and $\beta > 2$, then
\[ h(G) > 2, \quad h(F-LG) > 3 \quad \mbox{and} \quad  h(\b{F}) >5 \]
for every linear form $L(\x) \in \mathcal{O}_{K}[\x]$.
\end{lemma}

\begin{proof}
Let $\tau_{r}$ be the diagonal $n \times n$ matrix which has $\pi$ as its first $r$ entries and $1$ otherwise. If $h(G) \le 2$ then we may write 
\[ G = x_{1}L_{1}+x_{2}L_{2} \pmod \pi, \]
for some linear forms $L_{i} \in \mathcal{O}_{K}[\x]$. If we let 
\begin{displaymath}
T =  \left( \begin{array}{cc}
1 & 0 \\
0 & \pi^{-1} 
\end{array} \right),
\end{displaymath}
then $\b{F'} = T\b{F_{\tau_{2}}}$ is defined over $\mathcal{O}_{K}$. However $\beta - 2 > 0$, contradicting our assumption that $\b{F}$ is $(\alpha,\beta)$-reduced (cf. condition (\ref{bot-cond}), p.\pageref{bot-cond}). If $h(F-LG) \le 3$ for some linear form $L(\x) \in \mathcal{O}_{K}[\x]$, then we may write 
\[ F-LG = x_{1}Q_{1}+x_{2}Q_{2}+x_{3}Q_{3} \pmod \pi, \]
for some quadratic forms $Q_{i} \in \mathcal{O}_{K}[\x]$. If we let 
\begin{displaymath}
T =  \left( \begin{array}{cc}
\pi^{-1} & -\pi^{-1}L \\
0 & 1 
\end{array} \right),
\end{displaymath}
then $\b{F'} = T\b{F_{\tau_{3}}}$ is defined over $\mathcal{O}_{K}$. However $\alpha-3>0$, contradicting that $\b{F}$ is $(\alpha,\beta)$-reduced. Finally if $h(\b{F}) \le 5$ then we may write 
\[ \b{F'}=(F-LG,G) = (x_{1}Q_{1}+\ldots+x_{5}Q_{5},x_{1}L_{1} + \ldots + x_{5}L_{5}) \pmod \pi, \]
for some quadratic and linear forms $Q_{i},L_{i},L \in \mathcal{O}_{K}[\x]$. This time we let 
\begin{displaymath} 
T =  \left( \begin{array}{cc}
\pi^{-1} & 0 \\
0 & \pi^{-1} 
\end{array} \right),
\end{displaymath}
then $\b{F''}=T\b{F'}_{\tau_{5}}$ is defined over $\mathcal{O}_{K}$. However $\alpha+\beta-5>0$, contradicting that $\b{F}$ is $(\alpha,\beta)$-reduced.
\end{proof}


\section{Reduced systems}

In this section we will work with our system $\b{F}=(F,G)$ modulo $\pi$, which from now on we shall denote as $\b{f}=(f,g)$. We assume that $\b{F}$ is $(\alpha,\beta)$-reduced, for some $\alpha>3$ and $\beta>2$. Hence $\b{f}$ will satisfy the conclusion of Lemma \ref{reduced} viz.
\begin{equation} \label{h-invcon}
h(g) > 2 \quad h(f-lg) > 3 \quad \mbox{and} \quad  h(\b{f}) >5 
\end{equation}
for every linear form $l(\x) \in \F_{q}[\x]$. We also denote $m = o(\b{f})$.

The aim of this section is to show that we can find a non singular zero of $\b{f}$ which by Hensel's Lemma will lift to give us a zero of our original system $\b{F}$. For clarity we outline the steps we will take in order to prove this: 
\\[0.05cm]

\b{Step 1:} We prove that we can find a zero $\b{e}_{1}$ say, of $\b{f}$ such that $\nabla g \ne 0$. 

Therefore we are able to write our system $\b{f}$ in the shape
\begin{eqnarray*}
f(\x) & = & x_{1}^{2}f_{1} + x_{1}f_{2} + f_{3} \\
g(\x) & = & x_{1}g_{1} + g_{2}
\end{eqnarray*}
where $g_{1} \not\equiv 0$. If $\b{e}_{1}$ is a non singular zero then we're done. Otherwise we can find some $\lambda \in \F_{q}$ such that $f_{1} = \lambda g_{1}$. We now consider the equivalent system: 
\begin{eqnarray*}
(f-\lambda x_{1} g)(\x) & = & x_{1}(f_{2} - \lambda g_{2}) + f_{3} \\
g(\x) & = & x_{1}g_{1} + g_{2}.
\end{eqnarray*}
We may therefore assume that $\b{f}$ is equivalent to one of two situations:
\begin{itemize}
\item[(i)] $\deg_{x_{1}} f = 0$,
\item[(ii)] $\deg_{x_{1}} f = 1$.
\end{itemize}

\b{Step 2:} In case (i) we show that we may write $\b{f}$ as 
\begin{eqnarray*}
f(\x) & = & f(\x_{2}) \\
g(\x) & = & x_{1}x_{2} + g_{2}(\x_{3})
\end{eqnarray*}
where we define $\x_{i} := (x_{i},x_{i+1},\ldots,x_{m})$. Next we will find a non singular zero $\x_{2}$ of $f$ such that $x_{2} \ne 0$. Hence by setting $x_{1} = x_{2}^{-1}g_{3}(\x_{3})$, we get a non singular zero of the system $\b{f}$ as required.

\b{Step 3:} In case (ii) we show that we can write $\b{f}$ as 
\begin{eqnarray*}
f(\x) & = & x_{1}f_{2}(\x_{3}) + f_{3}(\x_{2}) \\
g(\x) & = & x_{1}x_{2} + g_{2}(\x_{3}).
\end{eqnarray*}
Next we define the quartic form
\[ H(\x_{2}) := x_{2}f_{3}(\x_{2}) - (f_{2}g_{2})(\x_{3}). \]
It will follow that if we can find a non singular zero of $H$ such that $x_{2} \ne 0$ then we can find a non singular zero of the system $\b{f}$. Finding a non singular zero of $H$ such that $x_{2} \ne 0$ requires a blend of ideas which utilizes the information (\ref{h-invcon}) we have at our disposal about the $h$-invariant of the system. 

Having described the outline of the proof we proceed with Step 1.


\subsection*{Step 1}

Throughout this step and subsequent steps we will need to make use of three important Lemmata, the first being attributed to Warning \cite{War35} (for example see \cite[Theorem 1E, p.137]{Sch76}).

\begin{lemma} \label{war35}
Let $F_{1},\ldots,F_{r}$ be a system of forms of degrees $d_{1},\ldots,d_{r}$ respectively in $m$ variables over $\F_{q}$. If $m > \delta = \sum d_{i}$, then the system $F_{1},\ldots,F_{r}$ has at least $q^{m-\delta}$ common affine $\F_{q}$-rational zeros.
\end{lemma}

The next Lemma comes from a book of Schmidt \cite[Lemma 3A, p.147]{Sch76}.

\begin{lemma} \label{schbook}
Let $F$ be a non zero polynomial over $\F_{q}$ in $m$ variables of total degree $d$. Then the number $N_{a}$ of affine zeros of $F$ in $\F_{q}^{n}$ satisfies 
\[ N_{a} \le dq^{m-1}. \] 
\end{lemma}

The final lemma in our toolbox is due to Leep \& Yeomans \cite{LeeYeo94}. 

\begin{lemma} \label{LeeYeo}
Let $P \in \F_{q}[x,y]$ be an absolutely irreducible polynomial of degree $d$. Then the number $N$ of non-singular zeros of $P$ satisfies
\begin{eqnarray*}
N \ge q + 1 - \frac{1}{2} (d-1)(d-2) [2\sqrt{q}], 
\end{eqnarray*}
where $[\gamma]$ denotes the least integer not exceeding $\gamma$.
\end{lemma}

\begin{proof}
Write $S$ to denote the number of $\F_{q}$ singular zeros of $P$. Then if the curve defined by $P(x,y)=0$ has genus $g$ (not to be confused with the quadratic form $g$), it follows from Corollary $1$ of Leep \& Yeomans \cite{LeeYeo94} that
\[ | N + S -(q+1)| \le g([2\sqrt{q}]-1) + \frac{1}{2}(d-1)(d-2). \]
Next we use the above estimate together with the following bound on the genus
\[ g \le \frac{1}{2}(d-1)(d-2)-S, \] 
which comes from Lemma 1 of \cite{LeeYeo94}, to obtain the required bound
\[ N \ge q+1 - \frac{1}{2}(d-1)(d-2)[2\sqrt{q}]. \]
\end{proof}

We say that the system $\b{f}$ is equivalent the system $\b{f'} = (f',g')$, if
\begin{eqnarray*}
f'(\x) &=& (f - lg)(\tau \x)  \\
g'(\x) &=& g(\tau \x)
\end{eqnarray*}
for some linear form $l \in \F_{q}[\x]$ and $\tau \in GL(n,\F_{q})$. That being said we let $a \ge 0$ denote the maximum integer such that we can write
\begin{eqnarray*}
f(\x) &=& f(y_{1},\ldots,y_{a},z_{1},\ldots,z_{b}) \\
g(\x) &=& g(z_{1},\ldots,z_{b})
\end{eqnarray*}
where $o(\b{f})=m= a+b$ and $o(g)=b$, amongst all systems equivalent to $\b{f} = (f,g)$. We introduce the notation $\y = (y_{1},\ldots,y_{a})$ and $\z = (z_{1},\ldots,z_{b})$. Throughout this step we shall assume that every zero of $\b{f}$ is such that $\z = \b{0}$, otherwise the zero will be such that $\nabla g \ne \b{0}$, which completes Step 1.
Since $m \ge h(\b{f}) \ge 6$, then by Lemma \ref{war35} there exists a zero $\b{e}_{1}$ say of $\b{f}$. Also note that since $\b{f}(x_{1},0,\ldots,0) = 0$ for all $x_{1} \in \F_{q}$, then $m-1 \ge h(\b{f}) \ge 6$. Let  $N(\b{f})$ denote the number of affine zeros of the system $\b{f}$ over $\F_{q}$. Hence by Lemma \ref{war35},
\[ N(\b{f}) \ge q^{2}. \]
So there exists another zero of $\b{f}$ not in the affine span of $\b{e}_{1}$, say $\b{e}_{2}$.
If all the zeros of $\b{f}$ are in the affine span of $\{ \b{e}_{1},\b{e}_{2} \}$ and $q > 2$ then $m-2 \ge h(\b{f}) \ge 6$. Hence by Lemma \ref{war35}, 
\[ N(\b{f}) \ge q^{3}. \]
Otherwise there is a zero of $\b{f}$ not in the span of $\{ \b{e}_{1},\b{e}_{2} \}$. In either case we may assume that there are at least $3$ linearly independent zeros of $\b{f}$, say $\b{e}_{1},\b{e}_{2},\b{e}_{3}$.

If any of these zeros are such that $\nabla g \ne \b{0}$ then we have completed Step 1. Otherwise we may assume each of these zeros $\b{e}_{i}$ are singular for $g$, which implies that $a \ge 3$. Note that we may write
\[ f(\x) = f_{a}(\y) + f_{a,b}(\y,\z) \]
where $f_{a}(\y) :=f(\y,\b{0})$. 

Suppose that $b \ge 5$, then by Lemma \ref{war35}
\[ N(\b{f}) \ge q^{a+b-5} \ge q^{a}. \]
If all the zeros of $\b{f}$ are such that $\z = \b{0}$, then by the above inequality $\b{f}$ must vanish in the span of $\{ \b{e}_{1},\ldots,\b{e}_{a} \}$. Hence we have $m - a \ge h(\b{f}) \ge 6$, so that 
\[ N(\b{f}) \ge q^{a+1} \]
Therefore we conclude that in any case we can find a zero of $\b{f}$ such that $\z \ne \b{0}$, completing Step 1 provided $b \ge 5$. Moreover since $b \ge h(g) \ge 3$ we can find a zero $\b{e}_{a+1}$ say of $g$, therefore as before we have that $b-1 \ge h(g) \ge 3$. Hence we may assume from now on that $b = 4$, otherwise as the above argument shows $b \ge 5$ would allow us to find a zero with $\z \ne \b{0}$.

Next we would like to find a zero $\b{e}_{1}$ say of $f_{a}$, such that the variable $x_{1}$ appears in $f_{a}$.
Let $a' = o(f_{a})$ and write $f_{a}(\y) = f_{a}(\y')$, where $\y' = (y_{1},\ldots,y_{a'})$ after a non singular change of the variables $\y$. Note the following relationship
\[ 6 \le h(\b{f}) \le h(f_{a}) + b = h(f_{a}) + 4. \]
Therefore $a' \ge h(f_{a}) \ge 2$. We now show that we can find a zero of $\b{f}$ such that $\y' \ne \b{0}$. Suppose all zeros of $\b{f}$ were such that $\y' = \b{0}$ and $\z = \b{0}$, that leaves $a-a'$ variables which are non zero for each solution of $\b{f} = \b{0}$. However
\[ N(\b{f}) \ge q^{a' + (a-a') + b -5} \ge q^{a-a' +1}. \]
Which implies that either we can find a zero of $\b{f}$ such that $\z \ne \b{0}$ (completing Step 1) or we can find a zero such that $\y' \ne \b{0}$. Hence we may assume that $\b{e}_{1}$ is a non trivial zero of $f_{a}$. So if $\b{e}_{a+1}$ is a zero of $g$ we can assume that $f(\b{e}_{a+1}) \ne 0$ otherwise we would have found the required zero to complete this step. So we have that
\[ f(\x) = y_{1}f_{2}(y_{1},\ldots,y_{a}) + f_{3}(y_{2},\ldots,y_{a}) + f_{a,b}(\y,\z), \]
where $f_{2} \not\equiv 0$ and $f_{a,b}(\b{e}_{a+1}) \ne 0$. By Lemma \ref{schbook} for $q>2$, we can find a vector $\b{e}_{2}$ say such that $f_{2}(\b{e}_{2}) \ne 0$. Finally we consider the following slice of the cubic form.
\begin{equation} 
S(X,Y,Z) := f(X\b{e}_{1} + Y\b{e}_{2} + Z\b{e}_{a+1}) = X^{2}u(Y,Z) + Xv(Y,Z) + w(Y,Z), 
\end{equation}
where $v(Y,0) \ne 0$ and $w(0,Z) \ne 0$. Note that if we can find a zero of $S$ such that $Z \ne 0$, then we can have found a zero of $\b{f}$ with $\z \ne \b{0}$ which completes Step 1. Note the following Lemma.

\begin{lemma} \label{lemS}
$S(X,Y,Z)$ has a zero with $Z \ne 0$.
\end{lemma}

\begin{proof}
If $S$ is absolutely irreducible we can set $Z=1$ and apply Lemma \ref{LeeYeo}, to deduce the existence of a zero over $\F_{q}$ for all $q$. If $S$ is reducible then we may assume that either it is the product of $3$ conjugate linear factors or the product of a linear factor defined over $\F_{q}$ and a quadratic factor. 

It cannot be the product of $3$ conjugate linear factors otherwise a $X^{3}$ term would appear in $S$. Hence we may assume that $S$ is the product of a linear factor over $\F_{q}$ and a quadratic factor. Note that $Z$ does not divide $S$, since $v(Y,0) \ne 0$. Therefore we can choose $(X,Y,Z)$ so that we set the linear factor equal to $0$ and have $Z \ne 0$, completing the proof.
\end{proof}

This completes Step 1 viz. we can find a vector $\b{e}_{1}$ such that $\b{f} = \b{0}$ and $\nabla g \ne \b{0}$.


\subsection*{Step 2}

In this step we may assume that our zero $\b{e}_{1}$ of $\b{f}$ is such that $\nabla g \ne \b{0}$ and $\deg_{x_{1}} f = 0$. Therefore we may immediately write
\begin{eqnarray*}
f(\x) & = & f(\x_{2}) \\
g(\x) & = & x_{1}g_{1}(\x_{2}) + g_{2}(\x_{2}),
\end{eqnarray*}
where $g_{1} \not\equiv 0$. By a non singular change of variables we can assume that $g_{1} = x_{2}$. Hence we may write
\begin{eqnarray*}
g(\x) & = & x_{1}x_{2} + g_{2}(\x_{2}) \\
& = & x_{2}(x_{1} + \lambda x_{2} + L(\x_{3})) + \hat{g}_{2}(\x_{3}),
\end{eqnarray*}
for some constant $\lambda \in \F_{q}$ and linear form $L(\x_{3})$. So by mapping $x_{1}$ to $x_{1} - \lambda x_{2} - L(\x_{3})$ and writing $g_{2}$ to denote $\hat{g}_{2}$ we have
\begin{eqnarray*}
f(\x) & = & f(\x_{2}) \\
g(\x) & = & x_{1}x_{2} + g_{2}(\x_{3}).
\end{eqnarray*}

Our goal now is to find a non-singular zero of $f(\x_{2})$ such that $x_{2} \ne 0$. Then by setting $x_{1} = -x_{2}^{-1}g_{2}(\x_{3})$, we obtain a non-singular zero of our system $\b{f}$, as required. We can appeal to a result of Lewis \& Schuur \cite[Theorem 3]{LewScr73} to immediately answer this exact question viz.

\begin{theorem}[Lewis \& Schuur, 1973] \label{lewscr73}
Let $k$ be a finite field of cardinality $q \ge 5$. Let $F$ be a non-degenerate cubic form over $k$ such that $o(F) \ge 4$ and let $L$ be a linear form over $k$. Then $F$ has a $k$-point which is a non-singular zero of $F$ and is not a zero of $L$.
\end{theorem}

To keep with the slicing theme we will prove our own version of Theorem \ref{lewscr73} by employing the information we have about the $h$-invariant. First we need a lemma.

\begin{lemma} \label{nscubic}
Let $F(x_{1},\ldots,x_{m})$ be a non-degenerate cubic form over  any finite field, such that $o(F) =m  \ge 4$. Then $F$ has a non-singular zero.
\end{lemma}

\begin{proof}
By Lemma \ref{war35}, $F$ has a non-trivial zero $\b{e}_{1}$ say. Therefore we may write
\[ F(\x) = x_{1}^{2}F_{1}(\x_{2}) + x_{1}F_{2}(\x_{2}) + F_{3}(\x_{2}). \]
If $F_{1} \not\equiv 0$ then $\b{e}_{1}$ is a non-singular zero. Otherwise,
\[ F(\x) = x_{1}F_{2}(\x_{2}) + F_{3}(\x_{2}) \]
where $F_{2} \not\equiv 0$ since $F$ is non-degenerate. In any finite field we can find a point $\x_{2} \in \F_{q}^{m-1}$ such that $F_{2}(\x_{2}) \ne 0$. Therefore we obtain a non-singular zero of $F$ by setting $x_{1} = (-F_{2}^{-1}F_{3})(\x_{2})$, as required.
\end{proof}

Next we need to consider $\delta := \deg_{x_{2}} f(\x_{2})$. If $\delta = 0$, then by Lemma \ref{nscubic} we can find a non-singular zero of $f$ and set $x_{2}=1$ and $x_{1} = -g_{2}$ to obtain non-singular zero of $\b{f}$. If $\delta = 1$ then
\[ f(\x_{2}) = x_{2}f_{2}(\x_{3}) + f_{3}(\x_{3}), \]
where $f_{2} \not\equiv 0$ since $\delta =1$ and $f_{3} \not\equiv 0$ since $h(f) > 3$. By Lemma \ref{schbook} we can find a vector $\x_{3} \in \F_{q}^{n-2}$ such that $f_{2}(\x_{3}),f_{3}(\x_{3}) \ne 0$, provided $q>5$. Therefore by setting $x_{2} = (-f_{2}^{-1}f_{3})(\x_{3}) \ne 0$ and $x_{1} = -x_{2}^{-1}g_{2}(\x_{3})$, we get a non-singular zero of $\b{f}$ as required. If $\delta = 2$, then $\b{e}_{2}$ is a non-singular zero of $\b{f}$ as required. Therefore we may assume that $\delta =3$.

By Lemma \ref{nscubic} we can find a non-singular zero $\b{a} = (a_{2},\ldots,a_{n})$ say of $f(\x_{2})$. If $a_{2} \ne 0$, then we can find a non-singular zero of $\b{f}$ by setting $x_{1} = -a_{2}^{-1}g_{2}(a_{3},\ldots,a_{n})$. So if $a_{2} = 0$, we can make a change of variables so that $\b{a} = \b{e}_{3}$. Therefore
\[ f(\x_{2}) = x_{3}^{2}f_{1}'(x_{2},\x_{4}) + x_{3}f_{2}'(x_{2},\x_{4}) + f_{3}'(x_{2},\x_{4}) + f_{3}''(\x_{4}), \]
where $\deg_{x_{2}} f_{3}' =3$ (since $\delta=3$), $f_{1}' \not\equiv 0$ (since $\b{e}_{3}$ is a non-singular zero) and $f_{3}'' \not\equiv 0$ (since $h(f) > 3$). By Lemma \ref{schbook} we can find a vector $\b{e}_{4}$ such that, $f_{2}'(x_{2},\b{e}_{4}), f_{3}''(\b{e}_{4}) \ne 0$, provided $q > 5$. We now consider the slice
\begin{eqnarray*}
T(X,Y,Z) &=& f(X\b{e}_{2} + Y\b{e}_{3} + Z\b{e}_{4}) \\
&=& cX^{3} + Y^{2}u(X,Z) + Yv(X,Z) + w(X,Z) + dZ^{3},
\end{eqnarray*}
where $c,d \ne 0$, $u(X,1) \ne 0$ and $\deg_{X} w \le 2$.

\begin{lemma}
$T(X,Y,Z)$ has a non-singular zero such that $X \ne 0$, provided $q>3$.
\end{lemma}

\begin{proof}
If $T$ is absolutely irreducible we can set $X=1$ and apply Lemma \ref{LeeYeo} to deduce the existence of a zero over $\F_{q}$ for all $q$. If $T$ is reducible then either it is the product of $3$ conjugate linear factors or the product of a linear factor defined over $\F_{q}$ and a quadratic factor. 

It cannot be the product of $3$ conjugate linear factors otherwise a $Y^{3}$ term would appear in $T$. Hence we may assume that $T$ is the product of a linear factor over $\F_{q}$ and a quadratic factor. Note that $X$ does not divide $T$, since $d \ne 0$. Therefore we can assume that the linear factor is $X - l(Y,Z)$, where $l \not\equiv 0$. Let $Q(X,Y,Z)$ denote the quadratic factor of $T$. If $X-l(Y,Z)$ does not divide $Q$ then $Q(l(Y,Z),Y,Z) \not\equiv 0$. By Lemma \ref{schbook} provided $q > 3$ we can find some $(Y,Z)$ such that $l(Y,Z), Q(l(Y,Z),Y,Z) \ne 0$. Therefore by setting $X = l(X,Y) \ne 0$ we get the required non-singular zero. On the other hand if $X-l(Y,Z)$  does divide $Q$ then we may write 
\[ T(X,Y,Z) = (X-l(Y,Z))^{2}(cX - l'(Y,Z)). \]
Now note that since $\deg_{Y} T = 2$, $cX - l'(Y,Z)$ cannot divide $X-l(Y,Z)$. So by letting $Q(X,Y,Z) =(X-l(Y,Z))^{2}$ we may find the required non-singular zero as before.
\end{proof}

This completes Step 2 viz. If $\b{e}_{1}$ is a zero of $\b{f}$ such that $\nabla g \ne \b{0}$ and $\deg_{x_{1}} f = 0$, then we can find a non-singular zero of $\b{f}$ provided $q > 5$.


\subsection*{Step 3}

In this step we shall assume that $\b{e}_{1}$ is a zero of $\b{f}$ such that $\nabla g \ne \b{0}$ and $\deg_{x_{1}} f = 1$. Therefore we can write
\begin{eqnarray*}
f(\x) &=& x_{1}f_{2}(\x_{2}) + f_{3}(\x_{2}) \\
g(\x) &=& x_{1}g_{1}(\x_{2}) + g_{2}(\x_{2})
\end{eqnarray*}
Recall (cf. beginning of Step 2) that we can make a change of basis so that 
\[ g(\x) = x_{1}x_{2} + g_{2}(\x_{3}), \]
and $\deg_{x_{1}} f =1$. Moreover by subtracting linear multiples of $g$ from $f$ we may assume that $f_{2} = f_{2}(\x_{3})$. As mentioned in the outline of this step, we define
\[ H(\x_{2}) = x_{2}f_{3}(\x_{2}) - (f_{2}g_{2})(\x_{3}). \]
Our strategy is to find a vector $\x_{2}$ to be able to apply the following lemma.

\begin{lemma} \label{nsH}
If we can find a non-singular zero $\x_{2} \in \F_{q}^{n-1}$ say of $H$ such that $x_{2} \ne 0$, then $(-x_{2}^{-1}g_{2}(\x_{3}),\x_{2}) \in \F_{q}^{n}$ is a non-singular zero of $\b{f}$.
\end{lemma}

\begin{proof}
First we show that $(-x_{2}^{-1}g_{2}(\x_{3}),\x_{2})$ is a zero of $\b{f}$. It is clear that $\x$ is a zero of $g$ since, $x_{1} = -x_{2}^{-1}g_{2}(\x_{3})$. Also since $x_{2} \ne 0$ we have
\begin{eqnarray*}
H(\x_{2}) &=& x_{2}(f_{3}(\x_{3}) - x_{2}^{-1}(g_{2}f_{2})(\x_{3})) \\
&=& x_{2}f(\x).
\end{eqnarray*}
Therefore $f(\x) = 0$. Next suppose $\x$ is a singular zero of $\b{f}$, then the following gradient vectors must be linearly dependent 
\begin{eqnarray*}
\nabla f(\x) &=& (f_{2}, f_{32}, x_{1}f_{2i} + f_{3i}) \\
\nabla g(\x) &=& (x_{2},x_{1},g_{2i})
\end{eqnarray*}
where $f_{2i},f_{3i},g_{2i}$ denotes $\frac{\partial f_{2}}{\partial x_{i}},\frac{\partial f_{3}}{\partial x_{i}},\frac{\partial g_{2}}{\partial x_{i}}$ respectively for $2 \le i \le n$. Consequently we have the following vector identity
\[ f_{2} \nabla g = x_{2} \nabla f. \]
Looking at the components of this identity we have that
\begin{eqnarray}
x_{1}f_{2}  &=& x_{2}f_{32} \label{gradid1} \\
f_{2}g_{2i} &=& x_{2}(x_{1}f_{2i} + f_{3i}) \label{gradid2}
\end{eqnarray}
for $3 \le i \le n$. Next we note the gradient vector of $H$,
\[ \nabla H(\x_{2}) = (f_{3} + x_{2}f_{32}, x_{2}f_{3i}-g_{2i}f_{2}-g_{2}f_{2i}). \]
Let $(\nabla H)_{j}$ denote the $j$th component of the vector $\nabla H$ for $1 \le j \le n-1$. Then by (\ref{gradid1}), 
\[ (\nabla H)_{1} = f_{3} + x_{2}f_{32} = f_{3} + x_{1}f_{2} = f =0. \]
Also by (\ref{gradid2}) for $2 \le j \le n-1$ we have,
\begin{eqnarray*}
(\nabla H)_{j} &=&  x_{2}f_{3i}-g_{2i}f_{2}-g_{2}f_{2i} \\
&=& x_{2}f_{3i} -g_{2}f_{2i} - x_{2}(x_{1}f_{2i} + f_{3i}) \\
&=& -f_{2i}(x_{1}x_{2}+g_{2}) = -f_{2i}g = 0.
\end{eqnarray*}
Therefore $\nabla H = \b{0}$, a contradiction. Hence $\x$ is a non-singular zero of $\b{f}$ as required.
\end{proof}

We shall now show that we can find a non-singular zero $\x_{2}$ say of $H$ such that $x_{2} \ne 0$. Of course from the outset it may be possible that $H$ is the product of a quadratic form $Q$ say with itself. Then if $H = 0$ we must have $Q = 0$, therefore
\[ \frac{\partial H}{\partial x_{i}} = 2Q\frac{\partial Q}{\partial x_{i}} = 0, \quad \mbox{for all} \; \; 2 \le i \le n \]
implying that every zero of $H$ is singular. We shall show that this cannot happen, more precisely we will prove the following.

\begin{lemma} \label{Habsirr}
Suppose the $h$-invariant condition $(\ref{h-invcon})$ viz.
\[ h(g) > 2 \quad h(f-lg) > 3 \quad \mbox{and} \quad  h(\b{f}) >5 \]
then the form
\[ H(\x_{2}) = x_{2}f_{3}(\x_{2}) - (g_{2}f_{2})(\x_{3}) \]
is necessarily absolutely irreducible over $\F_{q}$.
\end{lemma}

\begin{proof}
As an outline we will distinguish between the cases in which $H$ is either the product of two quadratic forms or the product of a linear form and absolutely irreducible cubic form. Suppose $H$ factors over $\bar{\F}_{q}$, therefore we can write
\begin{eqnarray*}
H(\x_{2}) &=& x_{2}f_{3}(\x_{2}) - (g_{2}f_{2})(\x_{3}) \\
&=& (A(\x_{3}) + x_{2}B(\x_{2}))(A'(\x_{3}) + x_{2}B'(\x_{2})) \\
&=& AA' + x_{2}(AB' + A'B) + x_{2}^{2}BB' 
\end{eqnarray*}
for some forms $A,A',B,B' \in  \bar{\F}_{q}[\x_{2}]$. Hence 
\begin{equation} \label{AA'}
-g_{2}f_{2}=AA' . 
\end{equation}
Also by looking at the $x_{2}$ coefficient we deduce 
\begin{equation} \label{f3}
f_{3} = AB' + A'B + x_{2}BB'. 
\end{equation}

\b{Case (a):} Suppose $H$ is the product of a linear and absolutely irreducible cubic factor over some extension $K$ say of $\F_{q}$. Then by considering the action of Gal$(K:\F_{q})$ on the factors, it is easy to see that each factor must be defined over $\F_{q}$. Suppose $(A+ x_{2}B)$ is the linear factor, then (\ref{AA'}) implies that $A$ divides either $g_{2}$ or $f_{2}$. It cannot divide $g_{2}$ since if it did then by setting the two linear forms $A,x_{1}=0$ we would have that $g = 0$. Hence $h(g) \le 2$, contradicting (\ref{h-invcon}). So $f_{2} = AL$ for some linear form $L$ defined over $\F_{q}$. Therefore $A' = -g_{2}L$. So by (\ref{f3}) we have
\begin{eqnarray*}
f(\x) &=& x_{1}f_{2} + f_{3} \\
&=& x_{1}AL + AB' - g_{2}LB + x_{2}BB'. 
\end{eqnarray*}
Therefore by setting the linear forms $A,L,x_{2}=0$ we deduce that $f=0$. Hence $h(f) \le 3$, contradicting (\ref{h-invcon}).

\b{Case (b):} Suppose $A = g_{2}$ and $A' = -f_{2}$. By (\ref{f3}) we can write
\begin{eqnarray*}
f(\x) &=& x_{1}f_{2} + f_{3} \\
&=& -x_{1}A' + AB' + A'B + x_{2}BB' \\
&=& B'g + (-x_{1}+B)(A'+x_{2}B')
\end{eqnarray*}
recalling that
\begin{eqnarray*}
g(\x) &=& x_{1}x_{2} + g_{2} \\ 
&=& x_{1}x_{2} + A.
\end{eqnarray*}
Since $A= g_{2}$ and $A' = -f_{2}$ are defined over $\F_{q}$, then either the factors $(A+x_{2}B)$ and $(A'+x_{2}B')$ are defined over $\F_{q}$ or the quadratic extension of $\F_{q}$. In the former case $B$ and $B'$ have coefficients in $\F_{q}$ and by setting $x_{1},B,B' =0$ we deduce that $h(f) \le 3$, a contradiction to (\ref{h-invcon}). In the latter case we may assume that $B$ and $B'$ are conjugates of each other over the quadratic extension of $\F_{q}$. Hence we may write
\[ B(\x_{2}) = l_{1}(\x_{2}) + \alpha l_{2}(\x_{2}) \]
and
\[ B(\x_{2}) = \alpha l_{1}(\x_{2}) + l_{2}(\x_{2}) \]
for some linear forms $l_{1},l_{2} \in \F_{q}[\x_{2}]$ and $\alpha \in \bar{\F}_{q}$. So if we set $l_{1},l_{2},x_{1}=0$ we deduce that $h(f) \le 3$, a contradiction as before. This completes Case (b).

Before moving onto the final case we shall make a few remarks. If rank$(g_{2}) \ge 3$ and $H$ is the product of two quadratic factors then we can assume that we're in Case (b). This is because the rank condition implies that $g_{2}$ is absolutely irreducible and therefore the condition (\ref{AA'}) forces $\{ A,A' \} = \{ f_{2} , -g_{2} \}$, up to a scalar multiple in $\F_{q}$. We are therefore left to deal with the case in which rank$(g_{2}) \le 2$.

\b{Case (c):} Suppose rank$(g_{2}) \le 2$. So $g_{2}$ is reducible over $\bar{\F}_{q}$ and we may write 
\[ g_{2}(\x_{3}) = g_{2}(x_{3},x_{4}) = l_{1}(x_{3},x_{4})l_{2}(x_{3},x_{4}) \]
for some linear forms $l_{1},l_{2}$ defined over $\bar{\F}_{q}$. We can assume that $H$ is the product of two quadratic factors since the other possibility is dealt with in Case (a). Either $\{ A,A' \} = \{ f_{2} , -g_{2} \}$, up to a scalar multiple in which instance we are in Case (b), or we can assume that $l_{1}$ divides $A$ and $l_{2}$ divides $A'$. Hence $x_{3},x_{4}=0$ implies that $A,A'=0$. So as before we write 
\begin{eqnarray}
f(\x) &=& x_{1}f_{2} + AB' + A'B + x_{2}BB' \\
g(\x) &=& x_{1}x_{2} + l_{1}l_{2}.
\end{eqnarray}
So by setting the linear forms $x_{1},x_{2},x_{3},x_{4}=0$ we deduce that $\b{f} = \b{0}$. Hence $h(\b{f}) \le 4$, a contradiction to (\ref{h-invcon}). 

This completes the proof of the lemma.
\end{proof}

We are now in a position to find a non-singular zero $\x_{2}$ for $H$ where $x_{2} \ne 0$ which by Lemma \ref{nsH} will imply that there is a non-singular zero of $\b{F}$. To do this we employ a slicing approach, following an idea used by Wooley \cite{Woo08} for the case of degree $7$ and $11$ forms. We do this owing to the sharp bounds that are available for point counting on curves over $\F_{q}$ opposed to hypersurfaces. 

Before stating the next Lemma, we shall need to introduce some notation. Let $L$ be a field and consider a polynomial $f \in L[x_{0},x_{1},\ldots,x_{n}]$. When $\b{\xi} \in L^{3n+1}$, we write $f|_{\b{\xi}}=f|_{\b{\xi}}(X,Y)$ to denote the sliced polynomial
\[ f(\xi_{0}+X,\xi_{1}+\xi_{n+1}X+\xi_{2n+1}Y,\ldots,\xi_{n}+\xi_{2n}X+\xi_{3n}Y). \]

Next we shall note the following result of Cafure \& Matera \cite{CafMat06}.

\begin{lemma} \label{CafMat}
Let $f \in \F_{q}[x_{0},\ldots,x_{n}]$ be an absolutely irreducible polynomial of degree $d \ge 2$. Then the number of slices $\b{\xi} \in \F_{q}^{3n+1}$, for which the polynomial $f|_{\b{\xi}}$ is not absolutely irreducible, is at most $\frac{1}{2}(3d^{4}-4d^{3}+5d^{2})q^{3n}$.
\end{lemma}

\begin{proof}
This is Corollary $3.2$ of \cite{CafMat06}.
\end{proof}

Finally by Lemma \ref{Habsirr} and Lemma \ref{CafMat} if $q \ge 296$ there exists a slice $H|_{\b{\xi}}$ of 
\[ H(\x) = x_{2}f_{3}(\x_{2}) - (g_{2}f_{2})(\x_{3}) \] 
which is an absolutely irreducible curve. Moreover on this slice we must have that the $x_{2}$ component is not identically zero otherwise $H|_{\b{\xi}}$ would factor into the product $(g_{2}|_{\b{\xi}})(f_{2}|_{\b{\xi}})$. Finally by Lemma \ref{LeeYeo} (p.\pageref{LeeYeo}), taking $q > 293$ is more than sufficient to ensure the existence of a non-singular zero of $H|_{\b{\xi}}$ for which $x_{2} \ne 0$ and so of $H$ as required. This completes the proof of Theorem \ref{mainres}

\newpage

\bibliography{CubQuad}

\begin{thebibliography}{10}

\bibitem{Art65}
E.~Artin.
\newblock {\em The collected papers of \mbox{Emil Artin}}.
\newblock Addison-Wesley, London, 1965.

\bibitem{AxKoc65}
J.~Ax and S.~Kochen.
\newblock Diophantine problems over local fields. {I}.
\newblock {\em Amer. J. Math.}, 87:605--630, 1965.

\bibitem{BirLewMur62}
B.~J. Birch, D.~J. Lewis, and T.~G. Murphy.
\newblock Simultaneous quadratic forms.
\newblock {\em Amer. J. Math.}, 84:110--115, 1962.

\bibitem{Bra45}
R.~Brauer.
\newblock A note on systems of homogeneous algebraic equations.
\newblock {\em Bull. Amer. Math. Soc.}, 51:749--755, 1945.

\bibitem{CafMat06}
A.~Cafure and G.~Matera.
\newblock Improved explicit estimates on the number of solutions of equations
  over a finite field.
\newblock {\em Finite Fields Appl.}, 12(2):155--185, 2006.

\bibitem{Dem50}
V.~B. Dem'yanov.
\newblock On cubic forms in discretely normed fields.
\newblock {\em Doklady Akad. Nauk SSSR (N.S.)}, 74:889--891, 1950.

\bibitem{Dem56}
V.~B. Dem'yanov.
\newblock Pairs of quadratic forms over a complete field with discrete norm
  with a finite field of residue classes.
\newblock {\em Izv. Akad. Nauk SSSR. Ser. Mat.}, 20:307--324, 1956.

\bibitem{Has24}
H.~Hasse.
\newblock \mbox{Darstellbarkeit von Zahlen durch quadratische Formen in einem}
  \mbox{beliebigen algebraischen Zahlkorper}.
\newblock {\em J. Reine Angew. Math.}, 153:113--130, 1924.

\bibitem{LeeYeo94}
D.~B. Leep and C.~C. Yeomans.
\newblock The number of points on a singular curve over a finite field.
\newblock {\em Arch. Math. (Basel)}, 63(5):420--426, 1994.

\bibitem{Lew52}
D.~J. Lewis.
\newblock Cubic homogeneous polynomials over {$p$}-adic number fields.
\newblock {\em Ann. of Math. (2)}, 56:473--478, 1952.

\bibitem{LewScr73}
D.~J. Lewis and S.~E. Schuur.
\newblock Varieties of small degree over finite fields.
\newblock {\em J. Reine Angew. Math.}, 262/263:293--306, 1973.
\newblock Collection of articles dedicated to Helmut Hasse on his seventy-fifth
  birthday.

\bibitem{Sch76}
W.~M. Schmidt.
\newblock {\em Equations over finite fields. {A}n elementary approach}.
\newblock Lecture Notes in Mathematics, Vol. 536. Springer-Verlag, Berlin,
  1976.

\bibitem{Sch80}
W.~M. Schmidt.
\newblock Simultaneous {$p$}-adic zeros of quadratic forms.
\newblock {\em Monatsh. Math.}, 90(1):45--65, 1980.

\bibitem{Spr55}
T.~A. Springer.
\newblock Some properties of cubic forms over fields with a discrete valuation.
\newblock {\em Nederl. Akad. Wetensch. Proc. Ser. A. {\bf 58} = Indag. Math.},
  17:512--516, 1955.

\bibitem{Ter66}
G.~Terjanian.
\newblock Un contre-exemple \`a une conjecture d'{A}rtin.
\newblock {\em C. R. Acad. Sci. Paris S\'er. A-B}, 262:A612, 1966.

\bibitem{War35}
E.~Warning.
\newblock \mbox{Bemerkung zur vorstehenden Arbeit von Herrn Chevalley}.
\newblock {\em Abh. Math. Sem. Hamburg}, 11:76--83, 1935.

\bibitem{Woo08}
T.~D. Wooley.
\newblock Artin's conjecture for septic and unidecic forms.
\newblock {\em Acta Arith.}, 133(1):25--35, 2008.

\end{thebibliography}
\bibliographystyle{plain}

\end{document}